\documentclass[a4paper,12pt]{article}
\usepackage{amssymb,latexsym,theorem,hyperref}\usepackage{graphicx}
\makeatletter
\newcommand{\F}{\mathbb F}
\newcommand{\A}{\mathbb A}
\newcommand{\Z}{\mathbb Z}
\newcommand{\Pp}{\mathbb P}

\newcommand{\End}{\mathop{\operator@font{End}}\nolimits}
\newcommand{\Hom}{\mathop{\operator@font{Hom}}\nolimits}
\newcommand{\GL}{\mathop{\operator@font{GL}}\nolimits}

\newcommand{\Spec}{\mathop{\operator@font{Spec}}\nolimits}
\newcommand{\dlog}{\mathop{\operator@font{dlog}}\nolimits}
\newcommand{\divisor}{\mathop{\operator@font{div}}\nolimits}
\newcommand{\res}{\mathop{\operator@font{res}}\nolimits}
\newcommand{\red}{{\mathop{\operator@font{red}}\nolimits}}
\newcommand{\cEnd}{\mathop{{\mathcal End}}\nolimits}
\newcommand{\cHom}{\mathop{{\mathcal Hom}}\nolimits}
\newcommand{\dR}{{\mathop{\operator@font{dR}}\nolimits}}

\newcommand{\id}{\mathop{\operator@font{id}}\nolimits}
\newcommand{\cO}{\mathcal O}
\newcommand{\cH}{{\mathcal H}}
\newcommand{\cL}{{\mathcal L}}
\newcommand{\cM}{{\mathcal M}}
\newcommand{\cI}{{\mathcal I}}

\makeatother

\newcommand{\qed}{\unskip\nobreak\hfill\hbox{ $\Box$}}
\newtheorem{Proposition}[subsection]{Proposition}
\newtheorem{Theorem}[subsection]{Theorem}
\newtheorem{Lemma}[subsection]{Lemma}

\theorembodyfont{\normalfont}
\newtheorem{Remark}[subsection]{Remark}

\newtheorem{Example}[subsection]{Example}
\newtheorem{Exercise}[subsection]{Exercise}

\title{Frobenius Splittings}
\author{Wilberd van der Kallen}
\date{}
\begin{document}\sloppy

\maketitle


\section{Introduction}
Frobenius splittings were introduced by V. B. Mehta and A. Ramanathan in \cite{Mehta-Ramanathan} and refined further by
S. Ramanan and Ramanathan in \cite{Ramanan-Ramanathan}.
Frobenius splittings have proven to be a amazingly effective when they apply.
Proofs involving Frobenius splittings tend to be very efficient. Other methods usually require a much more detailed
knowledge of the object under study. For instance, while showing that the intersection
of one union of Schubert varieties with another
union of Schubert varieties is reduced, one does not need to know where that intersection is situated, 
let alone what it looks like exactly. 

Before getting to serious applications we slowly  introduce the main concepts.

\section{Frobenius splittings for algebras}\label{ring case}
Fix a prime $p>0$.
Let $A$ be commutative ring of characteristic $p$. So $A$ contains the field $\F_p$ with $p$ elements.
The Frobenius homomorphism $\phi:A\to A$ is the ring map sending $a$ to $a^p$. The same notation $\phi$ 
will be used for the Frobenius 
homomorphism on other $\F_p$-algebras.
Let $\F$ be a perfect field of characteristic $p$. So the Frobenius map $\phi:\F\to\F$ is a field automorphism.
The field $\F$ will serve as our base field.

\paragraph{Pull back}
If $M$ is an $A$-module, then $\phi^\ast M$ denotes the $A$-module obtained by \emph{base change} along $\phi$.
That is, as an abelian group  $\phi^\ast M$ equals $M$, but there is a different module structure, given as follows.
\label{pullback module}Let us use $\lozenge $ 
to denote the new module structure. One puts
$$a\lozenge m:=a^pm\quad \mbox{for }a\in A,\ m\in \phi^\ast M.$$ If we interpret $\phi:A\to A$ as a map
 $\phi:A\to \phi^*A$, then $\phi$ is $A$-linear:
$$\phi(a+b)=\phi(a)+\phi(b),\qquad \phi(ab)=a\lozenge \phi(b).$$ 

\paragraph{Splitting}
Define a \emph{Frobenius splitting} on $A$ to
be an $A$-linear  map $\sigma:\phi^*A\to A$ with $\sigma\circ \phi=\id$.
In other words, $\sigma$ is a left inverse of $\phi$, whence the name Frobenius splitting.
We often just say \emph{splitting}. When $A$ has a splitting $\sigma$ we call $A$ or $\Spec(A)$ \emph{split} by $\sigma$.

\ 

\noindent
 A Frobenius splitting $\sigma$ of $A$ is just a set map $\sigma$ from $A$ to itself, satisfying
\begin{enumerate}
 \item $\sigma(a+b)=\sigma(a)+\sigma(b)$, for $a,b\in A$,
\item $\sigma(a\lozenge b)=a \sigma(b)$, for $a,b\in A$,
\item $\sigma(1)=1$.
\end{enumerate}
Notice that these three properties do imply $\sigma(\phi(a))=a$, because $\sigma(\phi(a))=\sigma(a\lozenge 1)=a \sigma(1)=a$. 

Call a map $\sigma:A\to A $  a \emph{twisted} linear endomorphism if it satisfies  (1) and (2).
Write $\End_\phi(A)$ for the abelian group of twisted linear endomorphism of $A$. We make it into
an $A$-module by putting $(a\ast \sigma)(b)=\sigma(ab)$ for $a\in A$, $\sigma\in\End_\phi(A)$, $b\in A$.
So the module structure on $\End_\phi(A)$ is given by premultiplication.  Postmultiplication as in 
$(a \sigma)(b)=a\sigma(b)$ describes the $A$-module structure on 
$\phi^*\End_\phi(A)$. 

Here is the first result. Recall that a ring is \emph{reduced} if it has no nonzero 
\emph{nilpotent} elements \cite[p. 33]{Eisenbud book}.

\begin{Lemma}\label{reduced algebra}
 If $A$ has a Frobenius splitting, then $A$ is reduced.
\end{Lemma}
\paragraph{Proof}If not, there is an $a\in A$, $a\neq0$ with $a^2=0$. But then $a=\sigma(\phi(a))=\sigma(a^p)=\sigma(0)=0$.
Contradiction.\qed

\

We can see this lemma as a first indication that possessing a Frobenius splitting is something special. 
After all, not all $A$ are reduced.

\paragraph{Polynomial rings}
We wish to understand  $\End_\phi(A)$ when $A$ is a polynomial ring  $\F[x_1,\dots,x_n]$ over our perfect field $\F$. 
Let us start with
the one variable case $A=\F[x]$.
The $A$-module $\phi^*A$ has a basis $1, x,\dots, x^{p-1}$, so $\sigma\in \End_\phi(A)$ is determined by
the $\sigma(x^i)$ with $i=0,\dots, p-1$.
Define $\sigma_0\in \End_\phi(A)$ by stipulating that $\sigma_0(x^{p-1})=1$, $\sigma_0(x^i)=0$ for $0\leq i< p-1$.

\begin{Lemma}\label{free}
 $\End_\phi(\F[x])$ is free with basis $\sigma_0$.
\end{Lemma}
\paragraph{Proof}
Let $\sigma\in\End_\phi(\F[x])$. Put $f_i=\sigma(x^i)$ for $0\leq i\le p-1$. We claim that $\sigma=\sum_{i=0}^{p-1}(f_{i}\lozenge x^{p-1-i})\ast\sigma_0$. 
Indeed, for $0\leq j\leq p-1$ one gets
\begin{eqnarray*}
 \left(\sum_{i=0}^{p-1}(f_{i}\lozenge x^{p-1-i} )\ast\sigma_0\right)(x^j)&=&\\
\sum_{i=0}^{p-1}\sigma_0(f_{i}\lozenge x^{p-1-i} x^j)&=&\\
\sum_{i=0}^{p-1}f_{i}\sigma_0(x^{p-1-i}x^j)&=&f_j
\end{eqnarray*}
\qed

\paragraph{Tensor products}
Let $A$, $B$ be $\F$-algebras. As $\F$ is perfect, there is a natural map from $\End_\phi(A)\otimes_\F\End_\phi(B)$
to $\End_\phi(A\otimes_\F B)$. For $\sigma\in\End_\phi(A)$, $\tau\in\End_\phi(B)$, $a\in A$, $b\in B$, we put
$(\sigma\otimes\tau)(a\otimes b)=\sigma(a)\otimes\tau(b)$.  This defines a twisted endomorphism $\sigma\otimes\tau$
of $A\otimes_\F B$.

\begin{Exercise}\label{one generator}
If $A=\F[x_1,\dots,x_n]$ then $\End_\phi(A)$ is free with basis $\sigma_0$, where 
$\sigma_0(x_1^{p-1}\cdots x_n^{p-1})=1$, while $ \sigma_0(x_1^{m_1}\cdots x_n^{m_{n}})=0$ if at least one
 $m_i+1$ is not divisible by $p$.
\end{Exercise}

\begin{Exercise}\label{homogeneous}
The algebra  $A=\F[x_1,\dots,x_n]$ is graded with each $x_i$ having degree one.
The element $\sigma_0$ of the previous exercise sends homogeneous polynomials to homogeneous polynomials.
If $f\in A$ is homogeneous, then $f*\sigma_0$ also sends homogeneous polynomials to homogeneous polynomials.
In particular, if $(f*\sigma_0)(1)$ has constant term $1$ and $f$ is homogeneous, then $f*\sigma_0$ is a splitting.
\end{Exercise}

\begin{Lemma}\label{coefficients} The following are equivalent
\begin{itemize}
 \item $f\ast\sigma_0$ is a splitting,
\item The coefficient of $x_1^{p-1}\cdots x_n^{p-1}$ in $f$ is one, and the other monomials
$x_1^{m_1}\cdots x_n^{m_{n}}$
with nonzero coefficient in $f$ have at least one
 $m_i+1$  not divisible by $p$.
\end{itemize}\qed
\end{Lemma}
\begin{Remark}
 The coefficient of $x_1^{p-1}\cdots x_n^{p-1}$ in $f$  is the value of $(f*\sigma_0)(1)$ at the origin.
\end{Remark}

\paragraph{Compatible ideal}
Let $\sigma\in\End_\phi(A)$ and let $I$ be an ideal of $A$.
We say that $\sigma$ is \emph{compatible} with $I$ if $\sigma(I)\subset I$. 

Write $\End_\phi(A,I)$ for $\{\;\sigma\in\End_\phi(A)\mid \sigma(I)\subset I\;\}$.
Clearly, if $\sigma$ is compatible with $I$ it induces
a map $\bar\sigma:A/I \to A/I$ that also satisfies (1) and (2). So we get a map $\End_\phi(A,I)\to\End_\phi(A/I)$.
It sends splittings to splittings.

\begin{Lemma}\label{radical}
 If $A$ has a Frobenius splitting compatible with $I$, then $I$ is a radical ideal.
\end{Lemma}
\paragraph{Proof}Indeed, $A/I$ is reduced by Lemma \ref{reduced algebra}.\qed

\paragraph{Localization}
If $S$ is a multiplicatively closed subset of $A$, not containing zero, 
consider the \emph{localization} $S^{-1}A$  of $A$ \cite[2.1]{Eisenbud book}. 
Recall that an element of $S^{-1}A$ may be written in
more than one way as a fraction $a/b$. There is a natural localization map
 $\End_\phi(A)\to\End_\phi(S^{-1}A)$, say $\sigma\mapsto\sigma_S$, where $\sigma_S(a/b)= \sigma(ab^{p-1})/b$ 
for $a\in A$, $b\in S$. Check that $\sigma_S$ is well defined. 
The localization map sends splittings to splittings. If $S$ contains no zero divisors, then $A$ is a subring of $S^{-1}A$ and
$\sigma$ is the restriction of $\sigma_S$ to $A$.

\paragraph{Completion}
If the ideal $I$ is finitely generated, then one checks that
any $\sigma\in \End_\phi(A)$ is continuous for the $I$-adic topology,
also known as the \emph{Krull topology}  \cite[7.5]{Eisenbud book}.
If $\hat A$ denotes the $I$-adic completion we get a map $\End_\phi(A)\to \End_\phi(\hat A)$.
It sends splittings to splittings. 

\begin{Lemma}\label{ringres}
 Let $f\in A$ be  a non zero divisor. Then $$\End_\phi(A,(f))=f^{p-1}*\End_\phi(A).$$
\end{Lemma}
\paragraph{Proof}
On the one hand, if $\sigma\in \End_\phi(A)$, then $(f^{p-1}*\sigma)(fa)=\sigma(f*a)=f\sigma(a)$ for $a\in A$,
so that $f^{p-1}*\sigma\in \End_\phi(A,(f))$. On the other hand, if $\sigma\in \End_\phi(A,(f))$ define $\tau:A\to A$ by
$\tau(a)=\sigma(fa)/f$. One checks that $\tau\in\End_\phi(A)$ and that $f^{p-1}*\tau=\sigma$.\qed

\begin{Example} The cross is split.\\[.4em]\mbox{} \qquad\includegraphics[scale=.1]{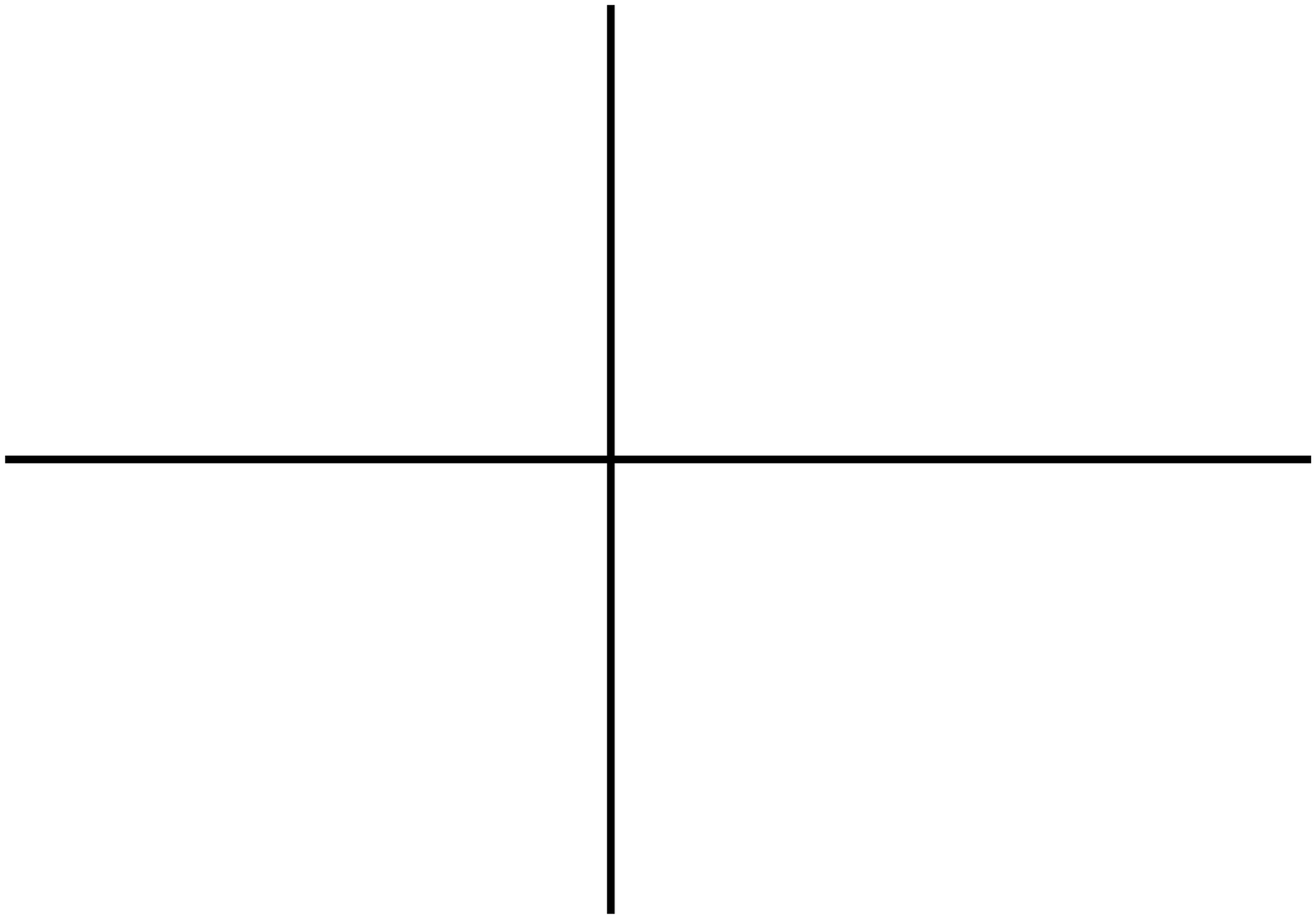}

 Let $A=\F[x,y]$ be the polynomial ring in two variables. The splitting $\sigma=(xy)^{p-1}\ast\sigma_0$ is compatible with the
ideal $(xy)$. Indeed, $\sigma(xyf)=\sigma_0(x^py^p f)=xy\sigma_0(f)$ for $f\in A$. So we have found a splitting on the
coordinate ring $\F[x,y]/(xy)$ of the union of the $x$-axis and the $y$-axis. This coordinate ring is not \emph{normal}
\cite[4.2]{Eisenbud book}.
The \emph{normalization} \cite[4.2]{Eisenbud book} is $\F[x]\times\F[y]$, and the map from the \emph{spectrum}
\cite[p. 54]{Eisenbud book} of $\F[x]\times\F[y]$ to the spectrum of $\F[x,y]/(xy)$
pinches together two points. So a Frobenius splitting does not rule out such behaviour. However, it does rule out 
pinching together two 
infinitely near points as displayed in the next example.
\end{Example}

\begin{Example} The cusp is not split.\\[.4em]\mbox{} \qquad\includegraphics[scale=.15]{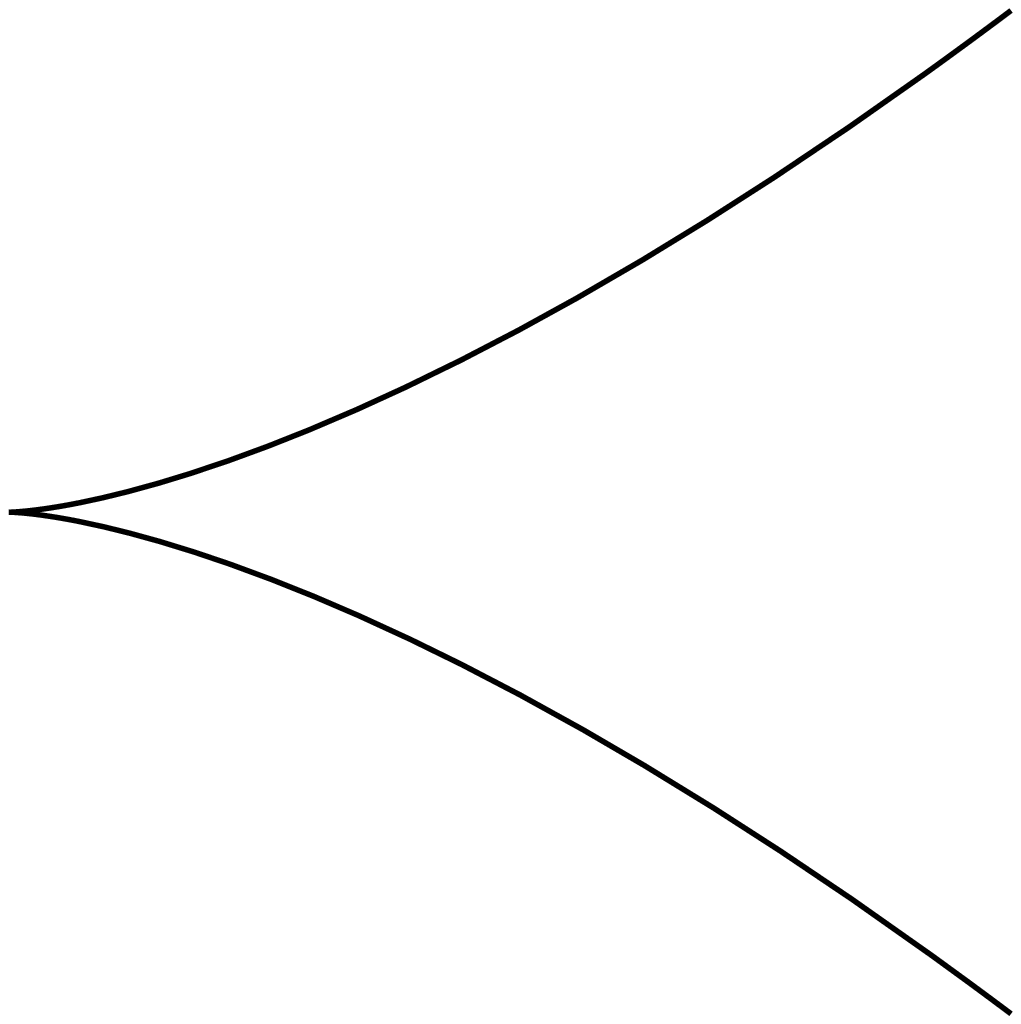}\label{cusp}

Consider the subring $A=\F[t^2,t^3]$ of the polynomial ring $\F[t]$. It is the coordinate ring of a cusp.
The polynomial ring $\F[t]$ is the normalization of $A$.
The ideal $\mathfrak c$ generated by $t^2$ and $t^3$ in $\F[t]$ is the \emph{conductor ideal} 
\cite[Exercise 11.16]{Eisenbud book}.
It is a common ideal in $A$ and in $\F[t]$. The ring $A/\mathfrak c$ is nonreduced. We already know that 
Frobenius splittings have little tolerance for nilpotents. So let us show that $A$ cannot have a splitting.
Suppose it did have a splitting $\sigma$. Take for $S$ the set of nonzero elements of $A$.
The splitting $\sigma_S$ on the field of  fractions $\F(t)$ must send $t^p$ to $t$. But it also should send $A$ to $A$.
Now $t^p$ is in $A$, but $t$ is not. Contradiction.
\end{Example}

\begin{Example} The node is split.\label{node}\\[.4em]\mbox{} \qquad \includegraphics[scale=.1,angle=0]{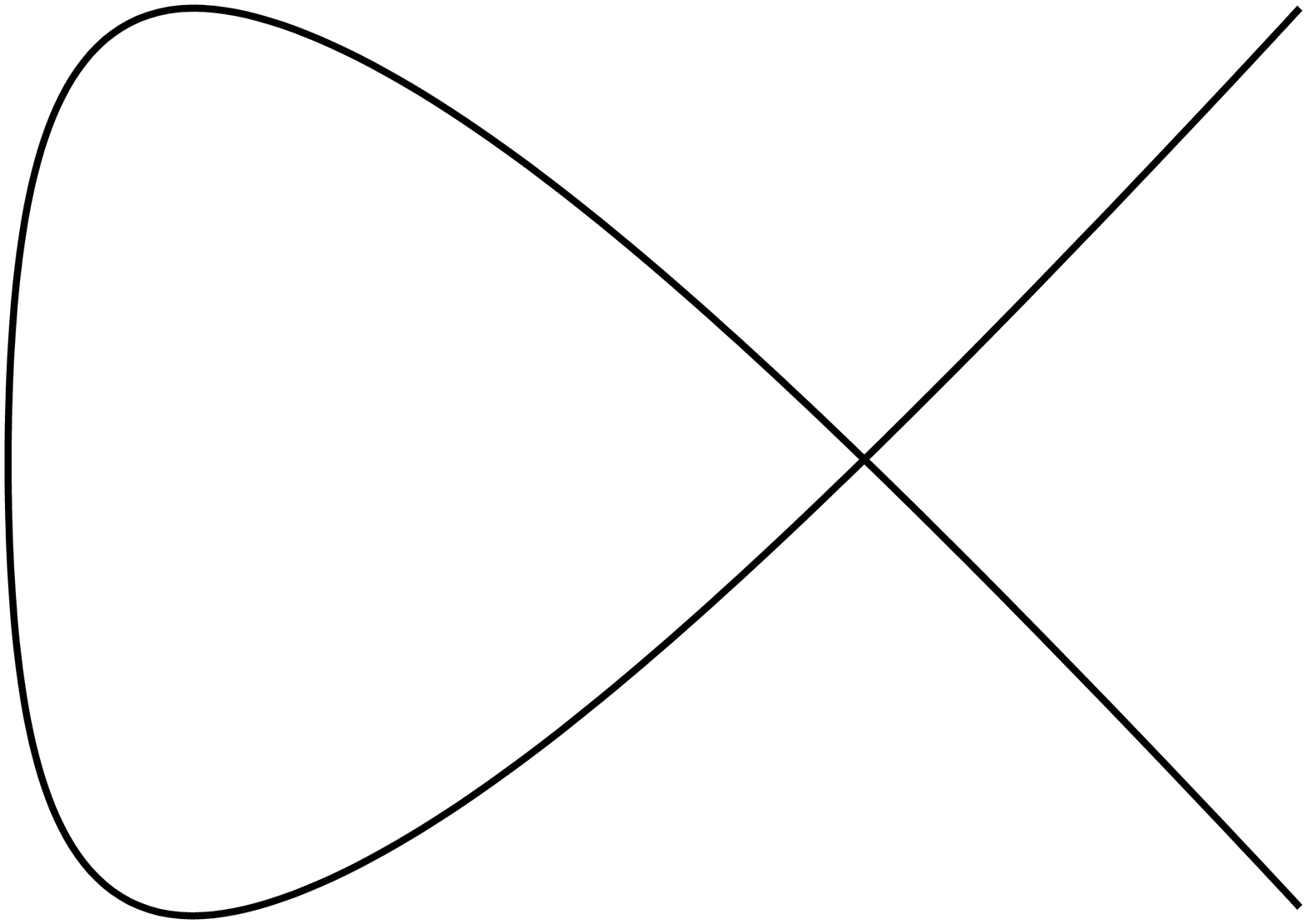}

Let our prime $p$ be unequal to two.
Consider the ring $A=\F[x,y]/(y^2-x^3-x^2)$, the coordinate ring of an ordinary node. 
(In characteristic two the equation $y^2=x^3+x^2$ would define a cusp.)

One may check by direct computation that $(y^2-x^3-x^2)^{p-1}\ast\sigma_0$ is  a splitting.
So the ring $A$ is Frobenius split compatibly with the ideal $(y^2-x^3-x^2)$.
For many purposes it is good enough to know just the existence of a splitting in $\End_\phi(\F[x,y],(y^2-x^3-x^2))$.
So then one would like to know that $1\in \F[x,y]$ is being hit by
the map from $\phi^*\End_\phi(\F[x,y],(y^2-x^3-x^2))$ to $\F[x,y]$ which sends
$\sigma$ to $\sigma(1)$. This is a linear algebra problem over a ring, so one has the usual tools of
localization and completion at one's disposal. But after localization and completion we see no difference between the node 
and the cross \cite[Second Example in 7.2]{Eisenbud book}. So this explains why the ideal of the node in the plane is compatibly split. 
\end{Example}

\begin{Example}No splitting when there is higher order contact.\\[.2em]\mbox{} \qquad 
\includegraphics[scale=.30,angle=0]{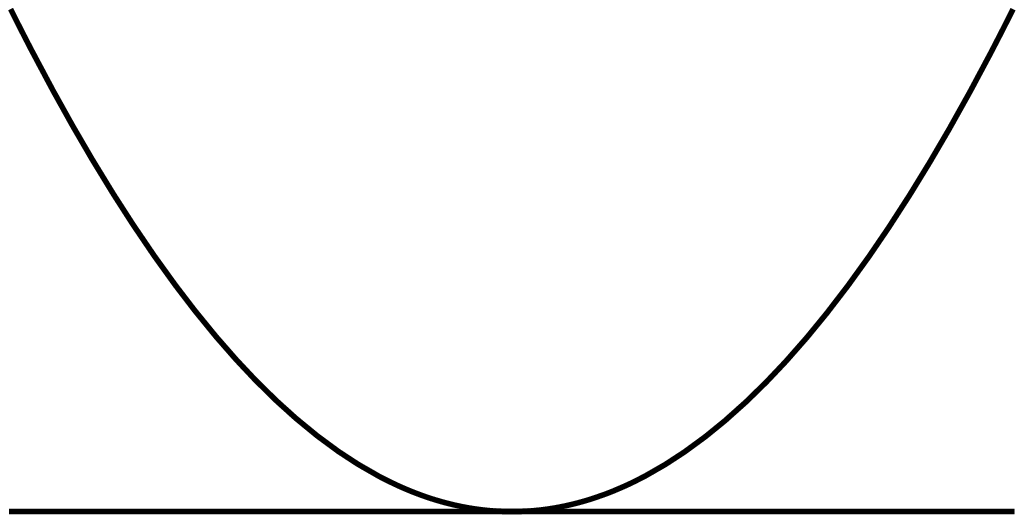}

We now look at the ideal $I=(y(y-x^2))\subset\F[x,y]$ of the union in the plane of the $x$-axis and the parabola $y=x^2$.
We claim this ideal is not compatibly split, the reason being the higher order contact at the intersection of
the parabola with the $x$-axis.
Suppose $\sigma\in\End_\phi(\F[x,y],I)$ were a splitting. 
First let $S$ consist of the powers of $y-x^2$. Then $\sigma_S\in \End_\phi(S^{-1}\F[x,y],S^{-1}I)$ maps 
$S^{-1}I$ to itself and also $\F[x,y]$ to itself. 

We claim the intersection of 
$\F[x,y]$ with  $S^{-1}I$ is the ideal $(y)$ of $\F[x,y]$. 
Indeed, a polynomial function on the plane that vanishes on an open dense subset of the $x$-axis vanishes on the whole
$x$-axis. Now take as open dense subset the intersection with the complement of the parabola.
So $\sigma$ is compatible with the ideal $(y)$ of $\F[x,y]$. Similarly, by inverting $y$ instead of $y-x^2$
we learn that $\sigma$ is compatible with the ideal $(y-x^2)$ in $\F[x,y]$ of
the other component. But then it must be compatible with the ideal $J=(y)+(y-x^2)$, the ideal of the scheme 
theoretic intersection
of the two components. However, because of the higher order contact, this scheme theoretic intersection is not reduced:
$\F[x,y]/J\cong \F[x]/(x^2)$ contains a nontrivial nilpotent. 
\end{Example}

\paragraph{Discussion}
What the last example shows us is that a Frobenius splitting allows to extrapolate from generic information,
on dense subsets of components, to information about a special locus. 
As V. B. Mehta explained it to me, a splitting seems to make bad behaviour at special points spread out to a neighborhood
of the bad point, which then makes it detectable generically. So behaviour that is not allowed generically
gets forbidden everywhere. In local coordinates one may think of the Frobenius map $t\mapsto t^p$ as a concentration flow, 
so that the splitting becomes a diffusion flow.

\section{Frobenius splittings for varieties}
For simplicity we take our base field $\F$ algebraically closed, still of characteristic $p$, $p>0$.
We will consider \emph{varieties} over $\F$, or more generally \emph{schemes} over $\F$ \cite{Hartshorne}. 
Unlike \cite{Brion-Kumar} or \cite{Hartshorne} we do not require
varieties to be smooth or connected. If $X$ is a variety over $\F$ we do require that $X$ is reduced and that
the corresponding morphism $X\to\Spec(\F)$ is of \emph{finite type} \cite[p. 84]{Hartshorne}.

\paragraph{Frobenius map for varieties}
So let $X$ be a variety or scheme over $\F$ with structure sheaf $\cO_X$. 
The absolute Frobenius map $F:X\to X$ is the morphism of ringed spaces
which is the identity on the underlying topological space and for any open subset $U$ raises a section 
$f\in \Gamma(U,\cO_X)$ to its $p$-th power.
So we use  $F$  as notation for a Frobenius map of schemes and $\phi$ for a Frobenius map of algebras. 
Note that $F$ is a map and that $\F$ is a field. Just like $\phi$ has been used for any Frobenius map of algebras, 
$F$ will be used for any Frobenius map of schemes.

\begin{Example}
 Let $A$ be an $\F$-algebra and $\phi$ its Frobenius endomorphism. The corresponding morphism $\Spec(A)\to \Spec(A)$ is
the absolute Frobenius map $F$.
\end{Example}

Say $X$ is a variety over $\F$. 
If one wants to view the morphism of ringed spaces $F:X\to X$ as a morphism of varieties over $\F$, 
then one may exploit the following commutative diagram of schemes
in which the vertical maps encode the $\F$-structure on $X$.

\[
\begin{array}{rcl}
{X} & \stackrel{F}{\longrightarrow} & {X} \\[.2em]
{\Big\downarrow} &
\mbox{
\Large$
\searrow$}

& {\Big\downarrow} 
\\
\Spec(\F)\hskip-1em & \stackrel{F}{\rightarrow} &\hskip-1em {\Spec(\F)} \\
\end{array}
\]

It suggests to view the source of $F:X\to X$ in a different way as a variety over $\F$, namely by using the diagonal
map instead of the vertical one.
The new variety structure obtained this way we denote $X^{(-1)}$.
Then $F:X^{(-1)}\to X$ is a map of varieties over $\F$. 

\paragraph{Splitting of a variety}
Let $X$ be a variety or scheme over $\F$. Consider the sheaf map $\cO_X\to F_*\cO_X$. 
Over any open subset $U$ it is described 
by the Frobenius map $\phi$
on the algebra $\Gamma(U,\cO_X)$. Following Mehta and \mbox{Ramanathan} we define a 
\emph{Frobenius splitting} $\sigma$ on $X$ to be a morphism of $\cO_X$-modules 
$F_*\cO_X\to\cO_X$ that splits the map $\cO_X\to F_*\cO_X$. 
So the composite $\cO_X\to F_*\cO_X\stackrel\sigma \to\cO_X$ must be the identity.
A scheme with a Frobenius splitting is called \emph{Frobenius split} or just \emph{split}.
Notice that the $\cO_X$-module structure on $F_*\cO_X$ is such that $\Gamma(U,F_*\cO_X)$ as a
$\Gamma(U,\cO_X)$-module is the pull back module
$\phi^*\Gamma(U,\cO_X)$ in the notation of Section
\ref{pullback module}. So a Frobenius splitting $\sigma$ on $X$ is a sheaf map that gives a Frobenius splitting of
 each algebra $\Gamma(U,\cO_X)$. Any open subset of a Frobenius split scheme is Frobenius split.

If $A$ is an $\F$-algebra, then $A$ is Frobenius split if and only if $\Spec(A)$
is Frobenius split. Thus Lemma \ref{reduced algebra} implies

\begin{Lemma}\label{reduced scheme}
 A Frobenius split scheme is reduced.\qed
\end{Lemma}
\begin{Exercise}\label{P1 splits}
 The splitting $x^{p-1}*\sigma_0$ of $\F[x]$ obtained from Lemma \ref{coefficients} gives a splitting
of $\A^1=\Spec(\F[x])$ that extends to a splitting of $\Pp^1$. 
\end{Exercise}

We write $\cEnd_F(X)$ for the sheaf of abelian groups
$\cHom(F_*\cO_X,\cO_X)$ associated to the presheaf $U\mapsto\End_\phi(\Gamma(U,\cO_X))$. 
We make it into an $\cO_X$-module using the $\Gamma(U,\cO_X)$-module structure $\ast$
on $\End_\phi(\Gamma(U,\cO_X))$ from Section \ref{ring case}. 
(Recall that we view $\sigma\in \End_\phi(\Gamma(U,\cO_X))$ as a map $\sigma: \Gamma(U,\cO_X)\to\Gamma(U,\cO_X)$ and
put $(a\ast \sigma)(b)=\sigma(ab)$.) 
A splitting is a global section $\sigma$ of $\cEnd_F(X)$ with $\sigma(1)=1$. If  $\sigma(1)$ is some nonzero constant
in $\F$, then we say that $\sigma$ \emph{spans a splitting}. The splitting it spans is of the form $\alpha*\sigma$
with $\alpha\in\F$.
Put $\End_F(X)=\Hom(F_*\cO_X,\cO_X)=\Gamma(X,\cEnd_F(X))$ and refer to its elements as twisted endomorphisms of $\cO_X$.

\begin{Example}
 If $A=\F[x_1,\dots, x_n]$, then $\Spec(A)=\A^n$ is Frobenius split. The sheaf $\cEnd_F(X)$ is a trivial line bundle
with nowhere vanishing
global section given by
the element $\sigma_0$ from Exercise \ref{one generator}. Locally any smooth $n$-dimensional variety looks like $\A^n$
up to completion, so $\cEnd_F(X)$ must be a line bundle on any smooth variety $X$. 
But there is no reason that it should be a trivial line
bundle and in fact it is usually not. We already see a problem with $\A^n$ itself. 
If $\alpha_1,\dots,\alpha_n\in\F$ are 
nonzero and $y_i=\alpha_ix_i$, 
then $\F[x_1,\dots, x_n]=\F[y_1,\dots, y_n]$. If $\sigma_0^x$ denotes the generator of $\End_\phi(A)$ constructed with the $x_i$
and $\sigma_0^y$  is the one constructed with the $y_i$, then $\sigma_0^x=(\alpha_1\cdots\alpha_n)^{p-1}\ast\sigma_0^y$.
That does not look like the transformation behaviour for structure sheafs.
\end{Example}

\begin{Theorem}[Mehta-Ramanathan]
 If $X$ is smooth of dimension $n$, then $\cEnd_F(X)$ is isomorphic as a line bundle to $\omega_X^{1-p}$, 
where $\omega_X$ is the canonical line bundle of $n$-forms on $X$ and $\omega_X^{1-p}$ means $\cHom( \omega_X^p,\omega_X)$.
\end{Theorem}
\begin{Remark}
 So if $\Gamma(X,\omega_X^{1-p})$ vanishes, then $X$ is certainly not split. Recall that a smooth complete variety is called 
a Fano variety if $\omega_X^{-1}$ is ample. While they are `quite rare' in algebraic geometry, they are very common in 
representation theory of reductive algebraic groups. 
Indeed, that is where many applications of Frobenius splittings are to be found.
\end{Remark}

The theorem is easy to prove with local duality theory, but we prefer to use the Cartier operator to construct a 
natural isomorphism between
the line bundles $\cEnd_F(X)$ and $\omega_X^{1-p}$. 
That will open the way to some explicit computations in local coordinates. This is important for checking
that a given global section of $\cEnd_F(X)$ is a splitting.
So let us digress and recall how the Cartier operator works.

\section{Cartier operator}
Let $X$ be a  variety of dimension $n$ over $\F$. We consider the DeRham complex
$$0\to\cO_X\to \Omega_X^1 \to \cdots \to \Omega_X^n \to 0
$$
with as differential $d$ the usual exterior differentiation.
Because this differential is not $\cO_X$-linear, we twist the $\cO_X$-module
structure on $\Omega_X^i$ by putting $f\lozenge \omega=f^p\omega$ for a section
$f\in \Gamma(U,\cO_X)$ and a differential $i$-form $\omega\in \Gamma(U,\Omega_X^i)$.
With this twisted module structure the DeRham complex is a complex of coherent
 $\cO_X$-modules, and the exterior algebra
$\Omega_X^*=\bigoplus_{i=0}^n\Omega_X^i$ is a differential graded
$\cO_X$-algebra.
We denote its cohomology sheafs $\cH_\dR^i$. They are $\cO_X$-modules by means of the twisted action.
If $U$ is an affine open subset, then $\Gamma(U,\cH_\dR^i)$ consists of all
closed differential $i$-forms on $U$ modulo the exact ones.
Now consider the map $\gamma: f\mapsto \hbox{class of }f^{p-1}df$ from
$\cO_X$ to $\cH^1_\dR$.

\begin{Lemma}\label{gamma}
$\gamma$ is a derivation and thus induces an $\cO_X$-algebra homomorphism
$c:\Omega_X^*\to \cH_\dR^*$.
\end{Lemma}
\begin{Remark}
Note that one should put the ordinary  $\cO_X$-module structure on
$\Omega_X^*$ here, not the twisted one that is used for $\cH_\dR^*$.
\end{Remark}

\paragraph{Proof of Lemma \ref{gamma}}
With
$$
\Phi(X,Y)=((X+Y)^p-X^p-Y^p)/p\in\Z[X,Y]
$$
we get
\begin{eqnarray*}
(f+g)^{p-1}d(f+g)&=&f^{p-1}df+g^{p-1}dg+d\Phi(f,g)\\
(fg)^{p-1}d(fg)&=&g\lozenge f^{p-1}df+f\lozenge g^{p-1}dg,
\end{eqnarray*}
where the first equality is a consequence of the fact that
$$
p(X+Y)^{p-1}d(X+Y)=pX^{p-1}dX+pY^{p-1}dY+pd\Phi(X,Y)
$$
in the torsion free $\Z$-module $\Omega^1_{\Z[X,Y]}$.
\qed

\begin{Proposition}\label{Cartier}
If $X$ is smooth, the homomorphism $c$ is bijective.

\end{Proposition}

\paragraph{Cartier operator}The inverse map $C:\cH_\dR^*\to\Omega^*_X$ is called the \emph{Cartier operator}
(cf. \rm\cite{Oesterle}).\index{Cartier operator}

\paragraph{Proof of Proposition \ref{Cartier}}
To check that a map of coherent sheafs is an isomorphism it suffices to
check that one gets an isomorphism after passing to the
completion at an arbitrary closed point. But then we are simply dealing with
the DeRham complex for a power series ring in $n$ variables over $k$ and
everything can be made very explicit (exercise).
\qed

\begin{Remark}
Here are some formulas satisfied by the Cartier operator, in informal
notation. In view of these
formulas the connection with Frobenius splittings is not surprising.
\begin{itemize}
\item $C(f^p\tau)=f C(\tau)$
\item $C(d\tau)=0$
\item $C(\dlog f)=\dlog f$, where $\dlog f$ stands for $(1/f) df$ if $f$
is invertible (or after $f$ has been inverted).
\item $C(\xi \wedge \tau)=C(\xi)\wedge C(\tau)$
\end{itemize}
Here $f$ is a function  and $\xi$, $\tau$ are forms.
\end{Remark}

\begin{Proposition}\label{Cartier duality}
If $X$ is smooth, we have a natural isomorphism of $\cO_X$-modules
$$\cEnd_F(X)\cong \omega_X^{1-p}=\cH om(\omega_X^p,\omega_X),$$
where $\omega_X$ is the canonical line bundle $\Omega_X^n$.
If $\tau$ is a local generator of $\omega_X$, $f$ a local section of
$\cO_X$, $\psi$ a local homomorphism $\omega_X^p\to\omega_X$,
then the local section $\sigma$ of $\cEnd_F(X)$ corresponding to $\psi$
is defined by $\sigma(f)\tau=C(\hbox{class of }\psi(f\tau^{\otimes p}))$.
\end{Proposition}
\paragraph{Proof}
One checks that $C(\hbox{class of }\psi(f\tau^{\otimes p}))/\tau$ does not
depend on the choice of $\tau$, so that $\sigma $ depends only on $\psi$.
To see that the map $\psi \mapsto \sigma$ defines an isomorphism of line
bundles we may argue as in the previous proof.
\qed

\begin{Remark}
 If $X$ is smooth of dimension zero, then $\omega_X=\cO_X$ and Proposition \ref{Cartier duality} describes 
the isomorphism $\End_\phi(\F)\cong\F$.
\end{Remark}

\begin{Example}
 We try the Proposition out for $X=\A^n=\Spec(\F[x_1,\dots,x_n])$. 
An obvious generator of $\omega_X$ is $\tau=dx_1\wedge\cdots\wedge dx_n$.
As local section of $\cO_X$ at some point we take a global function $f$ that does not vanish
at the point. The most obvious generator of $\Hom( \omega_X^p,\omega_X)$
sends $\tau^{\otimes p}$ to $\tau$. One may write it as $\tau^{1-p}$.
 We claim it corresponds with our old friend $\sigma_0$ from Exercise \ref{one generator}. We must check that 
$\sigma_0(f)\tau=C(f\tau)$, in simplified notation. It suffices to consider the case where $f$ is a monomial 
$x_1^{m_1}\cdots x_n^{m_{n}}$. 
If $m_n+1$  not divisible by $p$, then  $f\tau=dx_1\wedge\cdots\wedge dx_{n-1}\wedge d(x_nf)/(m_n+1)$ is a boundary so that
$C(f\tau)$ vanishes. Similarly $C(f\tau)$ vanishes if some other $m_i+1$ is not divisible by $p$.
So far the results are consistent with the definition of $\sigma_0$.
We still need to inspect the case where $x_1\cdots x_nf$ is a $p$-th a power, say $g^p$. Then $C(f\tau)$=
$C(g^p \dlog x_1\wedge\cdots \wedge\dlog x_n)=g\dlog x_1\wedge\cdots \wedge\dlog x_n=(x_1\cdots x_n)^{-1}g\tau$.
Indeed, $\sigma_0(f)$ equals $(x_1\cdots x_n)^{-1}g$.
\end{Example}

\section{Sheaf cohomology}

Let $X$ be a variety over $\F$ and let $\cL$ be a line bundle on $X$.
Then $F^*\cL$ is isomorphic to $\cL^p$. On an affine open subset $U$, say  $U=\Spec(A)$,
the isomorphism sends $a\otimes s\in A\otimes_{\phi}\Gamma(U,\cL)=\Gamma(U,F^*\cL)$ to 
$a(s^{\otimes p})\in \Gamma(U,\cL^p)$. 

If $\cM$ is a sheaf of abelian groups and $i\geq0$ then $H^i(X,F_*\cM)=H^i(X,\cM)$ as abelian groups, because $F$ is the
identity on the underlying topological space.

Now suppose that $\sigma$ splits $X$.
Then $\cL\otimes_{\cO_X}\cO_X\to \cL\otimes_{\cO_X}F_*\cO_X$ is split injective, so 
$H^i(X,\cL\otimes_{\cO_X}\cO_X)\to H^i(X,\cL\otimes_{\cO_X}F_*\cO_X)$ is split injective for each $i$.
By the projection formula $\cL\otimes_{\cO_X}F_*\cO_X$ equals $F_*(F^*\cL\otimes_{\cO_X}\cO_X)=F_*\cL^p$.
We get split injective maps $H^i(X,\cL)\to H^i(X,\cL^p)$. By iteration we get split
injective maps $H^i(X,\cL)\to H^i(X,\cL^{p^r})$
for $r\geq1$.

\def\citeMR1{\cite[Proposition 1]{Mehta-Ramanathan}}
\begin{Proposition}[\citeMR1]\label{MRprop1}
 Let $X$ be a projective variety which is Frobenius split. Let $\cL$ be a line bundle so that for some $i$,
$H^i(X,\cL^m)=0$ for  all large $m$ (e.g. $i>0$, $\cL$ ample). Then $H^i(X,\cL)=0$.\qed
\end{Proposition}

\def\citeMR2{Kodaira's vanishing theorem \cite[Proposition 2]{Mehta-Ramanathan}}
\begin{Proposition}[\citeMR2]
 Let $X$ be a smooth projective variety which is Frobenius split and $\cL$ an ample line bundle on $X$.
Then $H^i(X,\cL^{-1})=0$ for $i<\dim(X)$.
\end{Proposition}
\paragraph{Proof}

By Serre duality, $H^i(X,\cL^{-m})$ is the dual of $H^{n-i}(X,\omega\otimes\cL^{m})$ where $n=\dim(X)$.
Since $\cL$ is ample $H^{n-i}(X,\omega\otimes\cL^{m})$ vanishes for $i<\dim(X)$ and large $m$.
Thus $H^i(X,\cL^{-m})$ vanishes for $i<\dim(X)$ and large $m$. Now apply Proposition \ref{MRprop1}.\qed

\

We continue following \cite{Mehta-Ramanathan} for a while. Sometimes we are sketchy.
Let $Y\subset X$ be a closed subvariety.
We say that $\sigma\in\End_F(X) $ is \emph{compatible} with $Y$ if it maps the ideal sheaf $\cI_Y$ of $Y$ to itself.
Let $\End_F(X,Y)$ be the set of $\sigma\in \End_F(X)$ compatible with $Y$. 
If $X$ has a splitting compatible with $Y$ then we say that $Y$ is \emph{compatibly split}.
If a given splitting is compatible with several subvarieties, then we say that these subvarieties are
\emph{simultaneously compatibly split}.

\begin{Exercise}
 The splitting of $\Pp^1$ in Exercise \ref{P1 splits} is compatible with the points  $0$ and $\infty$. 
It corresponds with $\dlog(x)^{1-p}\in\Gamma(\Pp^1,\omega_{\Pp^1}^{1-p})$.
\end{Exercise}

\def\citeMR3{\cite[Proposition 3]{Mehta-Ramanathan}}
\begin{Proposition}[\citeMR3]\label{MR3}
 Let $X$ be a  projective variety and $Y\subset X$ a compatibly split closed subvariety.
If $\cL$ is an ample line bundle on $X$ then the restriction map $H^0(X,\cL)\to H^0(Y,\cL)$ is surjective and
$H^i(Y,\cL)$ vanishes for $i>0$.
\end{Proposition}
\paragraph{Sketch of proof}
Say $\sigma$ is the compatible splitting. Then $\sigma:F_*\cO_X\to \cO_X$ is split surjective and induces a split surjective
map $F_*\cI_Y\to \cI_Y$. Arguing as above we
 get a split surjective map from the 
$H^i(X,\cL^{p^r}\otimes_{\cO_X}\cI_Y)$ to the 
$H^i(X,\cL\otimes_{\cO_X}\cI_Y)$. Take $r$ large.
\qed

\begin{Lemma}
 Let $Y$ be a closed subvariety of $X$ and let $U$ be an open subset of $X$ such that $U\cap Y$ is dense in $Y$.
Then $\cEnd_F(X,Y)$ consists of the $\sigma\in \cEnd_F(X)$ whose restriction to $U$ lies in $\cEnd_F(U,U\cap Y)$.
In other words, compatibility with $Y$ may be checked on $U$.
\end{Lemma}
\paragraph{Proof} We have already used this principle in Example \ref{node}.\qed

\begin{Lemma}
 Let $Y$, $Z$ be closed subvarieties of $X$ and let $\sigma\in \cEnd_F(X)$.
If $\sigma$ is compatible with $Y$ and $Z$, then it is compatible with $Y\cup Z$
and with each irreducible component of $Y$. 
If $\sigma$ is compatible with $Y$ and $ Z$, then it is compatible with the scheme theoretic intersection 
 $Y\cap Z$.
If $Y$, $Z$ are simultaneously compatibly split, then their scheme theoretic intersection
is reduced.
\end{Lemma}
\paragraph{Hints for the Proof}Say $Y$ is irreducible.
By deleting the irreducible components of $Z$ that are not contained in $Y$ one forms an open subset $U$ of $X$
that intersects $Y$ in a dense subset and for which $U\cap (Y\cup Z)$ equals $U\cap Y$.
Let $V$ be open. A function $f\in\Gamma(V,\cO_X)$ vanishes on $V\cap (Y\cup Z)$ if and only if
it vanishes on both $V\cap Y$ and $V\cap Z$. 
The ideal sheaf of the scheme theoretic intersection $Y\cap Z$ is just $\cI_Y+\cI_Z$.\qed

\begin{Proposition}\label{three operations}
Let $X$ be a variety with Frobenius splitting $\sigma$.
 The collection of subvarieties with which $\sigma$ is compatible is closed under the following
operations 
\begin{itemize}
 \item Take irreducible components.
\item Take intersections.
\item Take unions.
\end{itemize}
\qed
\end{Proposition}

\def\citeMR4{\cite[Proposition 4]{Mehta-Ramanathan}}
\begin{Proposition}[\citeMR4]\label{push forward}
 Let $f:Z\to X$ be a proper morphism of algebraic varieties.
Assume that $f_*\cO_Z=\cO_X$. We then have
\begin{enumerate}
 \item If $Z$ is Frobenius split, then so is $X$.
\item If the closed subvariety $Y$ is compatibly split in $Z$, then so is $f(Y)$ in $X$.
\end{enumerate}

\end{Proposition}
\paragraph{Sketch of proof} The idea is that if $U$ is open in $X$, then a splitting of $\Gamma(U,\cO_X)$ amounts to the same as
a splitting of $\Gamma(f^{-1}(U),\cO_Z)$.\qed

\begin{Remark}\cite[Proposition 6.1.6]{bmod}
 The condition $f_*\cO_Z=\cO_X$ is satisfied when $f$ is surjective, proper, has connected fibres, and is separable in the 
following sense:
There is a dense subset of $x\in X$ for which there is $z\in f^{-1}(x)$ for which the tangent map $df_z$ is surjective,
where $df_z$ goes from the tangent space at $z$ to the tangent space at $x$.
A birational map is certainly separable.
\end{Remark}

\subsection{Residues.}\label{Residues}
Let $X=\Spec(A)$ be a smooth affine variety of dimension $n$ and let $f\in A$ have a smooth prime divisor of zeroes 
$D=\divisor(f)$.
There is a Poincar\'e residue map $\res:\Gamma(X,\omega_X(D))\to \Gamma(D,\omega_D)$.
Its composite with the surjective map $\beta\mapsto\beta\wedge\dlog(f):\Gamma(X,\Omega_X^{n-1})\to \Gamma(X,\omega_X(D))$ 
is 
the obvious restriction 
map $\Gamma(X,\Omega_X^{n-1})\to \Gamma(D,\omega_D)$. 
That characterizes the Poincar\'e residue.
We will not actually need the Poincar\'e residue, but we will
use the name \emph{residue} for some related maps, like
the map $f^{p-1}*\End_\phi(A)\to \End_\phi(A/(f))$ implied by Lemma \ref{ringres}. More specifically, 
we now seek an explicit formula for the $A$-linear residue map that is the composite of the maps\\
$f^{p-1}\Gamma(X,\omega_X^{1-p})\cong \End_\phi(A,(f))$,\\ $\End_\phi(A,(f))\to \End_\phi(A/(f))$,\\
$\End_\phi(A/(f))\cong \Gamma(D,\omega_D^{1-p})$.

\begin{Lemma}Take a point $P$ on $D$. Let $\tau$ be a local generator at $P$ of $\omega_D$ and let $\tilde \tau$ be 
a local lift to $\Omega^{n-1}_X$. For any sufficiently small neighborhood $U$ of $P$
 the $\Gamma(U,\cO_X)$-linear residue map $f^{p-1}\Gamma(U,\omega_X^{1-p})\to \Gamma(U\cap D,\omega_D^{1-p})$ sends 
$(\tilde\tau\wedge\dlog(f))^{1-p}$
to $\tau^{1-p}$.
\end{Lemma}
\paragraph{Proof}Take $U$ so small that $(\tilde\tau\wedge\dlog(f))^{1-p}\in f^{p-1}\Gamma(U,\omega_X^{1-p})$
and unravel the maps.\qed

\begin{Proposition}
 Let $X$ be smooth and let $\sigma\in\Gamma(X,\omega_X^{-1})$ so that $\sigma^{p-1}$ defines a splitting of $X$.
Then this splitting is compatible with the divisor $\divisor(\sigma)$ of $\sigma$.
\end{Proposition}
\paragraph{Proof}Take a smooth point $P$ on $\divisor(\sigma)_\red$, the reduced subscheme supporting $\divisor(\sigma)$. 
Locally around $P$ we are in the situation 
discussed above: $X=\Spec(A)$, $\sigma\in f^{p-1}\Gamma(X,\omega_X^{-1})$.\qed 

\begin{Remark}
 Actually $\divisor(\sigma)$ must be reduced. If there were a multiple component then one would get a vanishing residue there.
Or one could argue that $\divisor(\sigma)$ is a split scheme, hence reduced.
\end{Remark}

\begin{Lemma}
 Let $X$ be a smooth projective variety.
Let $\sigma\in\End_F(X)$. Suppose $\sigma(1)$ does not vanish at some point $P\in X$. Then 
 $\sigma$ spans a splitting of $X$. 
\end{Lemma}
\paragraph{Proof} As $\sigma(1)\in\Gamma(X,\cO_X)$, it is a constant function.\qed

\paragraph{Residually normal crossing}
The Lemma tells that if one has a section $\sigma$ of $\Gamma(X,\omega_X^{1-p})$, one may check if it spans a splitting 
by evaluating at a convenient point $P$. Now it happens often that there is a point 
$P$ where one may take a residue of $\sigma$ and thus bring the dimension down. Even better, one may have such luck that
by repeatedly taking residues the dimension can be brought all the way down to zero. And then finally,  in dimension zero,
one hopes to hit a nonzero constant.
That yields a rather practical way to establish that $\sigma$ spans a splitting.
The lucky situation we just alluded to has been formalized in \cite{MPL} with the notion `residually normal crossing'.

One may find it surprising that residually normal crossings are common. 
What happens is that in practice $\sigma$ is not chosen 
generically but in a very special position so as to have it compatible with an  effective divisor
that is important in the application at hand.
Then residually normal crossing may come as a bonus.
 
\begin{Example}
 We now give an example of the residual normal crossing phenomenon. It is not projective but affine. 
Actually one may extend the example to
nine dimensional projective space, but we will use Exercise \ref{homogeneous} instead. Let 
$X=\Spec(A)$ be the coordinate algebra of the space of $3$ by $3$ matrices. Thus $A=\F[x_{ij}]_{1\leq i,j\leq 3}$,
a polynomial ring in nine variables.
Take a volume form $\tau_9=d x_{11}\wedge dx_{12}\wedge \cdots \wedge dx_{33}$. The $9$ refers to the dimension.
Note that $\tau_9^{1-p}$ does not depend on how  we order the variables.
Now let $\sigma_9\in \Gamma(X,\omega_X^{1-p})$ be defined as\small{
$$\left((x_{11})\det\pmatrix{x_{11}&x_{12}\cr x_{21}&x_{22}}\det
\pmatrix{x_{11}&x_{12}&x_{13}\cr 
         x_{21}&x_{22}&x_{23}\cr
         x_{31}&x_{32}&x_{33}
}\det\pmatrix{ 
         x_{22}&x_{23}\cr
         x_{32}&x_{33}
}(x_{33})
\right)^{p-1}\tau_9^{1-p}.$$}
Taking a residue at the subvariety $x_{11}=0$ we get
$$\sigma_8=\left(\det\pmatrix{0&x_{12}\cr x_{21}&x_{22}}\det
\pmatrix{0&x_{12}&x_{13}\cr 
         x_{21}&x_{22}&x_{23}\cr
         x_{31}&x_{32}&x_{33}
}\det\pmatrix{ 
         x_{22}&x_{23}\cr
         x_{32}&x_{33}
}(x_{33})
\right)^{p-1}\tau_8^{1-p}.$$
We can take further residues in several ways, as the function in front of $\tau_8^{1-p}$ has several factors of the form
$x_{ij}^{p-1}$. Take residue at the subvariety $x_{12}=0$, then at its subvariety  $x_{21}=0$.
We arrive at $$\sigma_6=\left(\det
\pmatrix{0&0&x_{13}\cr 
         0&x_{22}&x_{23}\cr
         x_{31}&x_{32}&x_{33}
}\det\pmatrix{ 
         x_{22}&x_{23}\cr
         x_{32}&x_{33}
}(x_{33})
\right)^{p-1}\tau_6^{1-p}.$$
Take residue at $x_{22}=0$ and one is left with $\sigma_5=(x_{13}x_{31}x_{23}x_{32}x_{33})^{p-1}\tau_5^{1-p}.$
And so on until $\sigma_0=1$. What this means is that the original $\sigma_9$, viewed as an element of $\End_F(X)$
sends $1$ to a function $\sigma_9(1)$ with value $1$ at  the origin. But $\sigma_9$ is also given by a homogeneous
formula, so Exercise \ref{homogeneous} tells that $\sigma_9$ defines a splitting. 
It is compatible with the subvarieties that we encountered along the way. The splitting is also
compatible with all the other subvarieties that can be reached in a similar manner, such as the subvariety
$x_{11}=x_{33}=0$. But this is clear from Proposition \ref{three operations} anyway.
\end{Example}

\begin{Exercise}
 Extend the example to $n$ by $n$ matrices, or even $m$ by $n$ matrices.
\end{Exercise}

\section{Schubert varieties}
It is time to discuss a serious application. Mehta and Ramanathan constructed a Frobenius splitting on
the Bott-Samelson-Demazure-Hansen desingularisation of Schubert varieties in a flag variety $G/B$ to show that 
 the Schubert varieties are 
simultaneously compatibly split in the flag variety.
This then leads immediately to the result alluded to in the introduction about intersecting two unions of Schubert
varieties. This result about intersections is crucial in the analysis \cite{bmod}
of the fine structure (as $B$-modules) of 
dual Weyl modules, nowadays also known as costandard modules $\nabla(\lambda)$.
One also immediately gets that if $X$ is a Schubert variety in $G/B$, and $\cL$ is an ample line bundle,
then $\Gamma(G/B,\cL)\to \Gamma(X,\cL)$ is surjective and $H^i(G/B,\cL)=H^i(X,\cL)=0$ for $i>0$.
However, that is not quite the result that one wants. Kempf vanishing gives that in fact $H^i(G/B,\cL)=0$ for $i>0$ as
soon as $\Gamma(G/B,\cL)\neq0$. For instance, $H^i(G/B,\cO_{G/B})$ vanishes for $i>0$, but $\cO_{G/B}$ is not ample.
To get a result that covers all of Kempf vanishing, we will need the notion of $D$-splitting introduced by 
Ramanan and Ramanathan in \cite{Ramanan-Ramanathan}.

Let us first recall the Bott-Samelson-Demazure-Hansen resolution. 
We need the usual notations and terminology from the theory of reductive algebraic groups.
Let us remind the reader of some of the ingredients in a standard example. See also our book \cite{bmod} for more details
on these
constructions.

\begin{Example}Fix $n>1$.
 Let $G$ be the linear algebraic group $\GL_n$ over $\F$. As is common, we often discuss things as if $G$ is a group.
The true group is $G(\F)$, the group of $\F$-rational points of $G$. 
By $B$ we denote the algebraic subgroup of upper triangular matrices, by $T$ the algebraic
subgroup of diagonal matrices,
by $N(T)$ its normalizer, consisting of monomial matrices. (A monomial matrix is invertible and
has one nonzero entry in each row.)
The Weyl group $W=N(T)/T$ is isomorphic to the symmetric group on $n$ letters.
We let $S$ be the set of matrices that can be obtained by permuting two consecutive columns of the identity matrix.
One calls $S$ a set of representatives of fundamental reflections in the Weyl group. The number of elements of $S$ is $n-1$.
The algebraic subgroup of lower triangular matrices we denote $\tilde B$. So $B$, $\tilde B$ are opposite Borel
subgroups with intersection $T$.
The \emph{flag variety} $G/B$ parametrizes \emph{flags} in $n$-dimensional vector space.
Indeed, an invertible matrix $g$ defines a flag $L_1\subset L_2 \subset \cdots \subset L_n$ with $L_i$ being the span of
the first $i$ columns of $g$. Matrices $g$, $h$ define the same flag if and only if the cosets $gB$, $hB$ are equal.
For other reductive groups one still speaks of the flag variety $G/B$ in analogy with this example. 
For $s\in S$ let $P_s$ be the minimal
parabolic subgroup generated by $B$ and $s$. The subvariety $P_s/B$ of $G/B$ is isomorphic with a projective line $\Pp^1$.
A line bundle $\cL$ on $G/B$ is ample if and only if its restriction to $P_s/B$ is ample  for each $s\in S$.
There is a $G$-equivariant line bundle $\cL_\rho$ on $G/B$ so that $\cL_\rho^{-1}$ is `just ample', 
meaning that for each $s\in S$ its restriction to $P_s/B$ is the ample generator of the Picard group of $P_s/B$. 
One knows $\rho$ as the half sum of the positive roots.
A line bundle $\cL$ on $G/B$ is ample if and only
$\Gamma(G/B,\cL\otimes\cL_\rho)$ is nonzero. And $\Gamma(G/B,\cL)$ is nonzero if and only if the $\Gamma(P_s/B,\cL)$ are nonzero
for all $s\in S$.
\end{Example}

\paragraph{Bott-Samelson-Demazure-Hansen resolution}
If $X$, $Y$ are varieties with $B$ acting from the right on $X$ and from the left on $Y$, then the contracted product
$X\times^BY$ is the the quotient of $X\times Y$ by the equivalence relation $(xb,y)\approx (x,by)$ for $x\in X$, $y\in Y$,
$b\in B$, provided that quotient exists as a variety. 
If $s_1s_2\cdots s_d$ is a word on the alphabet $S$, so if one is given a sequence of length $d$ with values in
$S$, then we put $Z(s_1s_2\cdots s_d)=P_{s_1}\times^BP_{s_2}\times^B\cdots\times^BP_{s_d}/B$. 
Multiplication defines a map
from $Z(s_1s_2\cdots s_d)$ to $G/B$, sending $x_1\times^B\cdots\times^Bx_d B$ to $x_1\cdots x_d B$.
If the word is \emph{reduced} then $Z(s_1s_2\cdots s_d) \to G/B$ is birational to its image, which is of dimension $d$. 
This may be taken as a 
definition of \emph{reduced}. The image of $Z(s_1s_2\cdots s_d)$ in $G/B$ is the closure of a $B$-orbit. A $B$-orbit in 
$G/B$ is called a \emph{Schubert cell} and its closure is called a \emph{Schubert variety}. Schubert varieties may be singular.
The  closure of a $\tilde B$-orbit is called an \emph{opposite 
Schubert variety}.
The set of Schubert varieties is parametrized by $W$. If $\dot w\in N(T)$ is a representative of $w\in W$, 
then $X_w$ denotes the closure of the
orbit $B\dot wB/B$ of $\dot wB$. The dimension of $X_w$ is known as the length of $w$.
If $s_1s_2\cdots s_d$ is reduced, then it also describes an element $w\in W$
and the birational rational
map $Z(s_1s_2\cdots s_d)\to X_w$ is known as a Bott-Samelson-Demazure-Hansen resolution  of 
$X_w$. 
Indeed, $Z(s_1s_2\cdots s_d)$ is smooth: It is an iterated $\Pp^1$ fibration. The projection
map $Z(s_1s_2\cdots s_d)\to Z(s_1)=\Pp^1$, sending $x_1\times^B\cdots\times^Bx_d B$ to $x_1 B$
has fibre $Z(s_2\cdots s_d)$ above the point $B$ of $Z(s_1)$. We think of $B$ as point zero on this $\Pp^1$ and we think
of $s_1B$ as
the point $\infty$. On $Z(s_1s_2\cdots s_d)$ we have the divisor $Z_i$ consisting of the $x_1\times^B\cdots\times^Bx_d B$
with $x_i=1$. The divisors $Z_1$,\dots,$Z_d$ meet transversely at a point $P$. 
If $s_1s_2\cdots s_d$ is a reduced word of maximal length, then
 Mehta and Ramanathan show that $\cEnd_F(Z(s_1s_2\cdots s_d),Z_1\cup\cdots\cup Z_d)$ is the pullback from 
$G/B$ of $\cL_{\rho}^{1-p}$,
where $\cL_\rho^{-1}$ is the `just ample' line bundle on $G/B$. See also \cite[Proposition A.4.6]{bmod}, where the same is 
shown for any word, after Mathieu. 

\paragraph{A splitting}
The flag variety itself is also a Schubert variety. It corresponds with the longest element $w_0$ of the Weyl group.
Take a Bott-Samelson-Demazure-Hansen resolution
$Z(s_1s_2\cdots s_d)\to G/B$. (Although it is called a resolution, it is not a resolution of singularities, 
as $G/B$ itself is smooth.)
We wish to take a section $\tau\in\Gamma(G/B,\cL_\rho^{-1})$ which does not vanish in the image $B$ of the point $P$
where the $Z_i$ intersect each other.
A good choice for $\tau$ is a lowest weight vector, or simultaneous eigenvector for $\tilde B$,
 in $\Gamma(G/B,\cL_\rho^{-1})$. That works because the $\tilde B$-orbit of $B\in G/B$ is dense, so that $\tau$ cannot vanish
at the point $B$.
We will take $\tau$ this way.
Let $\sigma\in \cEnd_F(Z(s_1s_2\cdots s_d),Z_1\cup\cdots\cup Z_d)$ be the pullback of $\tau$. 
Then $\sigma^{p-1}$ spans a splitting because at $P$
there is a residually normal crossing. (A true normal crossing of the $Z_i$, actually.) 
This splitting is clearly compatible with the divisor $Z_1\cup\cdots\cup Z_d$.  
As $G/B$ is smooth, hence normal, the direct image 
of the structure sheaf of $Z(s_1s_2\cdots s_d)$ must be $\cO_{G/B}$, and Proposition \ref{push forward} gives a splitting
of $G/B$ compatible with the images of the $Z_i$. This covers all codimension one Schubert
varieties and with Proposition \ref{three operations} one shows that the splitting 
must be compatible with all Schubert varieties.

\begin{Theorem}[Mehta-Ramanathan]
 $G/B$ is Frobenius split with all Schubert varieties compatibly split.
\end{Theorem}

\begin{Remark}
Mehta and Ramanathan also considered Schubert varieties in $G/Q$ where $Q$ is a parabolic subgroup. 
\end{Remark}

\paragraph{Normality}We get a nice  proof of normality of Schubert varieties by means of the 
\begin{Lemma}[Mehta-Srinivas \cite{Mehta-Srinivas}]
 Let $f:Y\to X$ be a proper surjective morphism of irreducible $\F$-varieties.
Suppose that 
\begin{itemize}
\item $Y$ is normal,
\item the fibres of $f$ are connected,
\item $X$ is Frobenius split.
\end{itemize}
Then $X$ is normal.
\end{Lemma}
\paragraph{Discussion}The problem is local on $X$.
One argues as in Example \ref{cusp} (the example with the cusp) that if $f$ is in the function field of $X$
so that $f^p$ is a
regular function on some open $U$, then $f$ itself must be a regular function on $U$. 
That means that the map from the normalisation
of $X$ to $X$ cannot pinch together infinitely near points. In other words, one gets semi-normality in the sense of 
\cite{Swan}.
As the fibres of $f$ are connected, it is also impossible that disjoint points are pinched.
So $X$ is equal to its normalisation. In \cite[Proposition 1.2.5]{Brion-Kumar} the theme is 
worked out further by showing that every split scheme $X$
is weakly normal, meaning that every finite birational map $Z\to X$ is an isomorphism.

\

To apply the Lemma, one could show that a Bott-Samelson-Demazure-Hansen resolution
of a Schubert variety has connected fibres, but the argument in \cite{Mehta-Srinivas} is as follows. 
Let $X_w$ be a Schubert variety in $G/B$
and let $s_1s_2\cdots s_d$ be a corresponding reduced word.
Let $X_z$ be the image of $Z(s_1\cdots s_{d-1})$. By induction on dimension we may assume $X_z$ is normal.
With the Lemma one shows its image $X'$ in $G/P_{s_d}$ is normal. And the map from $X_w$ to $X'$ is a $\Pp^1$ fibration.
So $X_w$ is normal.

\begin{Theorem}
 Schubert varieties are normal.\qed
\end{Theorem}

\section{D-splittings}

To get more mileage out of the above construction of a splitting on $G/B$
one takes  a closer look at $\tau$ and $\cL_\rho^{-1}$.
We have not used yet that  $\cL_\rho^{-1}$ is ample. The line bundle $\cL_\rho^{-1}$ is well understood.
Recall that $\tau$ is a lowest weight vector in in $\Gamma(G/B,\cL_\rho^{-1})$. 
Its divisor $D$ is the union of the codimension one opposite Schubert varieties. (Compare \cite[Exercise 5.2.5]{bmod}.) 
Our splitting of $G/B$  is thus
simultaneously compatible with all Schubert varieties and all opposite Schubert varieties. 
But let us look at  cohomology.

\paragraph{D-splitting}
If $D$ is an effective divisor then a splitting $F_*\cO_X\to \cO_X$ of $X$ is called
a $D$-splitting    if it factors through
the map $F_*\cO_X \to F_*(\cO_X(D))$. So any $D$-splitting is a composite $F_*\cO_X \to F_*(\cO_X(D))\to \cO_X$.
If $X$ is smooth, and the section $\sigma$ of $\omega_X^{1-p}$ defines a splitting, then it is a
$D$-splitting precisely if $\sigma$ lands in
the subsheaf $ \omega_X^{1-p}(-D)$.
For example, in the above construction of the splitting on $G/B$ we may take for $D$
 the union of the codimension one opposite Schubert varieties.
If $X$ is $D$-split, then the surjective map $H^i(X,\cL^p)\to H^i(X,\cL)$ factors through $H^i(X,\cL^p\otimes \cO_X(D))$.
So if $i>0$ and $\cL^p\otimes \cO_X(D)$ is ample, then it factors through zero by Proposition \ref{MRprop1}.
We then conclude that $H^i(X,\cL)$ vanishes. Thus

\begin{Theorem}[Kempf vanishing]\label{Kempf vanishing}
Let $\cL$ be a line bundle on $G/B$ so that $\Gamma(G/B,\cL)$ is nonzero.
Then $H^i (G/B,\cL)$ vanishes for $i>0$.
\end{Theorem}
\paragraph{Proof}Indeed, with  $D$ as indicated above,
$\cL^p\otimes \cO_X(D)=\cL^p\otimes\cL_\rho^{-1}$ is ample.\qed

\

In  similar vein one wants to show

\begin{Theorem}\label{onto Schubert}
Let $\cL$ be a line bundle on $G/B$ so that $\Gamma(G/B,\cL)$ is nonzero.
Let $X_w$ be a Schubert variety in $G/B$.
Then $\Gamma(G/B,\cL)\to\Gamma(X_w,\cL)$ is surjective and $H^i (X_w,\cL)$ vanishes for $i>0$.
\end{Theorem}

\paragraph{Compatible D-splitting}
If $X$ is $D$-split and $Y$ is a subvariety of $X$ then we say that $Y$ is compatibly $D$-split if $Y$ is compatibly split 
and no irreducible component of $Y$ is contained in $D$. Assume this. The complement of $D$ intersects $Y$ in a dense 
open subset.

We claim that $F_*\cI_Y\to \cI_Y$  factors through
$F_*(\cI_Y(D))$. Indeed, $F_*(\cI_Y(D))\to \cO_X$ factors through $\cI_Y$,
 because a regular function on an open subset $U$ of $X$ vanishes on $U\cap Y$ if and only if 
it vanishes on a dense subset of $U\cap Y$. 

The surjective map $H^i(X,\cI_Y\otimes\cL^p)\to H^i(X,\cI_Y\otimes\cL)$ factors through 
$H^i(X,\cI_Y\otimes\cL^p\otimes \cO_X(D))$, and if $\cL^p\otimes \cO_X(D)$ is ample this vanishes for $i>0$,
by the proof of Proposition \ref{MR3}.

\paragraph{Proof of Theorem \ref{onto Schubert}} The Schubert variety is irreducible and contains the point $B\in G/B$ that 
lies in none of the opposite Schubert varieties. So we may argue as in the proof of
Theorem \ref{Kempf vanishing}. \qed

\section{Canonical splitting}The group $B$ acts on 
$\End_F(Z(s_1s_2\cdots s_d))=\Gamma(Z(s_1s_2\cdots s_d),\omega^{1-p})$ and one can check that 
our splitting is given by a $T$-invariant $\sigma$ in this $B$-module. Mathieu has observed that the $B$-module it generates
is rather small. So one might say the splitting is almost $B$-invariant. Mathieu has formalized this in the notion
\emph{canonical splitting} of a variety with $B$-action.

Recall that a $G$-module $M$ is called costandard if there is an equivariant line bundle $\cL$
on $G/B$ so that $M=\Gamma(G/B,\cL)$.
Mathieu employed canonical splittings to give an amazing proof of the following theorem
\begin{Theorem}
 The tensor product of two costandard modules
has a filtration by $G$-submodules whose associated graded module is a direct sum of costandard modules.
\end{Theorem}
See \cite{bmod}, \cite{Mathieu tilting} for an exposition of this.  

\section{More}There is much more that could be said, but we stop here. 
The Brion-Kumar book \cite{Brion-Kumar} is a treasure trove. If you want to see more recent
work, MathSciNet lists over forty references to \cite{Brion-Kumar}, and Google Scholar lists over a hundred.

\end{document}